\documentclass[12pt]{article}
\usepackage{amsfonts}

\usepackage{latexsym}
\usepackage{amssymb}
\usepackage{epsfig}
\usepackage{amsmath}
\usepackage{amsthm}
\usepackage[mathscr]{eucal}
\usepackage{subfigure}
\usepackage{graphicx}
\usepackage{hyperref}

\newtheorem{Lemma}{Lemma}

\newtheorem{Theorem}[Lemma]{Theorem}

\renewcommand{\qed}{\hfill{\ \ \rule{2mm}{2mm}} \vspace{0.2in}}

\begin{document}

\title{Existence of connected regular and nearly regular graphs}
\author{ \textbf{Ghurumuruhan Ganesan}
\thanks{E-Mail: \texttt{gganesan82@gmail.com} } \\
\ \\
New York University, Abu Dhabi }
\date{}
\maketitle

\begin{abstract}
For integers~\(k \geq 2\) and~\(n \geq k+1,\) we prove the following: If~\(n\cdot k\) is even, there is
a connected~\(k-\)regular graph on~\(n\) vertices. If~\(n\cdot k\) is odd, there is a connected nearly~\(k-\)regular
graph on~\(n\) vertices. 



\vspace{0.1in} \noindent \textbf{Key words:} Connected regular graphs, connected nearly regular graphs, existence.

\vspace{0.1in} \noindent \textbf{AMS 2000 Subject Classification:} Primary:
60J10, 60K35; Secondary: 60C05, 62E10, 90B15, 91D30.
\end{abstract}

\bigskip

\setcounter{equation}{0}
\renewcommand\theequation{\thesection.\arabic{equation}}
\section{Introduction}
The graph realization problem studies the existence of a graph on~\(n\) vertices
with a given degree sequence~\(d_1 \geq d_2 \geq \ldots \geq d_n.\)
Erd\H{o}s and Gallai (1960) provided necessary and sufficient conditions
for the sequence~\(\{d_i\}_{1 \leq i \leq n}\) to be graphic. Choudum (1986),
Aigner and Triesch (1994) and Tripathi and Vijay (2003) provide different
proofs of the Erd\H{o}s-Gallai theorem and Havel (1955) and Hakimi (1962) describe
algorithms for solving the graph realization problem.

If~\(d_i =k\) for all~\(1 \leq i \leq n,\) it is possible to use the Erd\H{o}s-Gallai Theorem
or the Havel-Hakimi criterion iteratively, to determine the existence of~\(k-\)regular graphs on~\(n\) vertices.
In this paper, we directly prove the existence of \emph{connected}~\(k-\)regular graphs
on~\(n\) vertices for all permissible values of~\(n \geq k+1,\) using induction.

We first have some definitions. Let~\(n \geq 2\) be any integer and let~\(K_n\)
the complete graph with vertex set~\(V = \{1,2,\ldots,n\}\)
and containing edges between all pairs of vertices. The edge between vertices~\(i\) and~\(j\) in~\(K_n\)
is referred to as~\((i,j).\) Let~\(G = (V,E)\) be any subgraph of~\(K_n\)
with edge set~\(E.\) Vertices~\(u\) and~\(v\) are \emph{adjacent} in~\(G\)
if the edge~\((u,v) \in E.\) A path~\(P = (v_1,\ldots,v_t)\) in~\(G\) is a sequence of distinct vertices such that~\(v_i\)
is adjacent to~\(v_{i+1}\) for~\(1 \leq i \leq t-1.\) The vertices~\(v_1\) and~\(v_t\)
are then connected by a path in~\(G.\) The graph~\(G\) is said to be \emph{connected} if any two vertices in~\(G\) are connected
by a path in~\(G\) (Bollobas~(2002)).

The degree of a vertex~\(v\) in the graph~\(G\) is the number of vertices adjacent to~\(v,\) in~\(G.\)
The graph~\(G\) is called a~\(k-\)\emph{regular} graph if the degree of each vertex in~\(G\)
is exactly~\(k.\) If~\(G\) is a~\(k-\)regular graph on~\(n\) vertices, then the product~\(n\cdot k\) is sum of degrees of the vertices in~\(G\)
which in turn is twice the number of edges in~\(G.\) Therefore it is necessary that either~\(n\) or~\(k\) is even. Finally, the graph~\(G\) is said to be a \emph{nearly}~\(k-\)\emph{regular} graph on~\(n\) vertices if~\(n-1\) vertices have degree~\(k\) and one vertex has degree~\(k-1.\) The following is the main result of the paper.
\begin{Theorem}\label{thm1} The following statements hold. \\
\((a)\) For all even integers~\(k \geq 2\) and all integers~\(n \geq k+1,\) there is a connected~\(k-\)regular graph on~\(n\) vertices.\\
\((b)\) For all odd integers~\(k \geq 3\) and all even integers~\(n \geq k+1,\) there is a connected~\(k-\)regular graph on~\(n\) vertices.\\
\((c)\) For all odd integers~\(k \geq 3\) and all odd integers~\(n \geq k+2,\) there is a connected nearly~\(k-\)regular graph on~\(n\) vertices.
\end{Theorem}

The paper is organized as follows. In Section~\ref{pf1}, we prove Theorem~\ref{thm1}.


\setcounter{equation}{0}
\renewcommand\theequation{\thesection.\arabic{equation}}
\section{Proof of Theorem~\ref{thm1}} \label{pf1}
The following fact is used throughout.\\
\((a1)\) If the minimum degree of a vertex in a graph~\(G\) is~\(\delta \geq 2,\)
then there exists a path in~\(G\) containing~\(\delta\) edges.\\
\emph{Proof of~\((a1)\)}: Let~\(P = (v_1,\ldots,v_t)\) be the longest path in~\(G;\) i.e.,~\(P\) is a path containing the maximum number of edges.
Since~\(P\) is the longest path, all the neighbours of~\(v_1\) in~\(G\) belong to~\(P.\) Since~\(v_1\)
has at least~\(\delta\) neighbours in~\(G,\) we must have~\(t \geq \delta+1\) and so~\(P\) has at least~\(\delta\) edges.~\(\qed\)

We prove~\((a),(b)\) and~\((c)\) in Theorem~\ref{thm1} separately below. \\
\emph{Proof of~\((a):\)} For~\(n = k+1,\) the complete graph~\(K_{n}\) is a connected~\(k-\)regular graph on~\(n\) vertices, with~\(k\) even. Suppose~\(G_{n}(k)\) is a connected~\(k-\)regular graph with vertex set~\(\{1,2,\ldots,n\},\) for some~\(n \geq k+1.\) Use property~\((a1)\) and let~\(P_n(k)\) be a path in~\(G_n(k)\) containing~\(k\) edges. Since~\(k\) is even, there are~\(\frac{k}{2}\) \emph{vertex disjoint} edges~\(e_i = (u_i,v_i), 1 \leq i \leq \frac{k}{2}\) in~\(P_n(k);\) i.e., the set~\(\bigcup_{1 \leq  i\leq \frac{k}{2}}\{u_i,v_i\}\) has~\(k\) distinct vertices. Remove the edges~\(e_i, 1 \leq i \leq \frac{k}{2}\) and add~\(k\) new edges~\[\bigcup_{i=1}^{\frac{k}{2}}\{(n+1,u_i)\}\bigcup \{(n+1,v_i)\}.\] The resulting graph is a connected~\(k-\)regular graph with vertex set\\\(\{1,2,\ldots,n+1\}\) and is defined to be~\(G_{n+1}(k).\)\\

\emph{Proof of~\((b):\)} The proof is analogous as in the case of~\((a).\) For~\(n = k+1,\) the complete graph~\(K_{n}\) is a connected~\(k-\)regular graph on~\(n\) vertices with~\(k\) odd. Suppose~\(G_{n}(k)\) is a connected~\(k-\)regular graph with vertex set~\(\{1,2,\ldots,n\},\) for some even~\(n \geq k+1.\) We use~\(G_{n}(k)\) to construct a connected~\(k-\)regular graph~\(G_{n+2}(k)\) with vertex set~\(\{1,2,\ldots,n+2\}\) as follows.

By property~\((a1),\) the graph~\(G_n(k)\) contains a path~\(Q_n(k) = (q_1,\ldots,q_{k+1})\) consisting of~\(k\) edges~\(\{(q_j,q_{j+1})\}_{1 \leq  j \leq k}.\) Remove the~\(k-1\) edges~\(\{(q_j,q_{j+1})\}_{1 \leq j \leq k-1}\) and add the following edges:\\
\((i)\) For~\(1 \leq j \leq k-1, j\) odd, add the edges~\(\{(n+1,q_j),(n+1,q_{j+1})\}.\)\\
\((ii)\) For~\(1 \leq j \leq k-1, j\) even, add the edges~\(\{(n+2,q_j),(n+2,q_{j+1})\}.\)\\
\((iii)\) Add the edge~\((n+1,n+2).\)

Since~\(k\) is odd, the total number of edges added in step~\((i)\) is~\(k-1\) and so there are~\(k-1\) edges with~\(n+1\) as an endvertex after step~\((i).\)
Similarly, after step~\((ii)\) there are~\(k-1\) edges with~\(n+2\) as an endvertex. Finally, the resulting graph after step~\((iii)\) is a connected~\(k-\)regular graph with vertex set~\(\{1,2,\ldots,n+1,n+2\}\) and is defined
to be~\(G_{n+2}(k).\)\\


\emph{Proof of~\((c):\)} Let~\(k\geq 2\) be odd and let~\(n \geq k+2\) be odd. We use the graph~\(G_{n-1}(k)\) with vertex set~\(\{1,2,\ldots,n-1\}\) obtained in~\((b)\) above to construct the graph~\(G_n(k)\) with vertex set~\(\{1,2,\ldots,n\}.\) From property~\((a1),\) the graph~\(G_{n-1}(k)\) contains a path~\(S_{n-1}(k)\) consisting of~\(k\) edges. Since~\(k\) is odd, there are~\(\frac{k-1}{2}\) vertex disjoint edges~\(f_i = (x_i,y_i), 1 \leq i \leq \frac{k-1}{2}\) in~\(S_{n-1}(k);\) i.e., the set~\(\bigcup_{1 \leq  i\leq \frac{k-1}{2}}\{x_i,y_i\}\) has~\(k-1\) distinct vertices. Remove the edges~\(f_i, 1 \leq i \leq \frac{k-1}{2}\) and add~\(k-1\) new edges~\[\bigcup_{i=1}^{\frac{k-1}{2}}\{(n,x_i)\}\bigcup \{(n,y_i)\}\] and define the resulting graph to be~\(G_{n}(k).\) By construction,~\(G_n(k)\) is connected, the vertex~\(n\) has degree~\(k-1\) and the rest of all the vertices have degree~\(k.\)~\(\qed\)



\subsection*{Acknowledgement}
I thank Professors Rahul Roy and Federico Camia for crucial comments and for my fellowships.

\bibliographystyle{plain}

\end{document}